\documentclass[10pt,reqno]{amsart}
\usepackage[hypertex]{hyperref}

\usepackage{amsmath,amsthm,amscd,amsfonts,amssymb}
\usepackage{euscript}

\numberwithin{equation}{subsection}

\newcommand{\sqsp}{\renewcommand{\baselinestretch}{1.1}\tiny\normalsize}

\newcommand{\relphantom[1]}{\mathrel{\phantom{#1}}}

\newtheorem{thm}[subsection]{Theorem}
\newtheorem{lemma}[subsection]{Lemma}

\newtheorem{corollary}[subsection]{Corollary}

\theoremstyle{definition}

\newtheorem{example}[subsection]{Example}

\newcommand{\bbC}{\mathbb{C}}

\newcommand{\UHLie}{U_{HLie}}

\DeclareMathOperator{\Id}{Id}
\DeclareMathOperator{\ad}{ad}
\DeclareMathOperator{\Hom}{Hom}

\begin{document}

\title{Hom-algebras and homology}
\author{Donald Yau}

\begin{abstract}
Classes of $G$-Hom-associative algebras are constructed as deformations of $G$-associative algebras along algebra endomorphisms.  As special cases, we obtain Hom-associative and Hom-Lie algebras as deformations of associative and Lie algebras, respectively, along algebra endomorphisms.  Chevalley-Eilenberg type homology for Hom-Lie algebras are also constructed.
\end{abstract}


\email{dyau@math.ohio-state.edu}
\address{Department of Mathematics, The Ohio State University at Newark, 1179 University Drive, Newark, OH 43055, USA}

\date{\today}
\maketitle
\sqsp


\section{Introduction}
\label{sec:intro}

In \cite{hls} Hartwig, Larsson, and Silvestrov introduced Hom-Lie algebras as part of a study of deformations of the Witt and the Virasoro algebras.  Closely related algebras also appeared earlier in the work of Liu \cite{liu} and Hu \cite{hu}.  A \emph{Hom-Lie algebra} is a triple $(L, \lbrack -,- \rbrack, \alpha)$, in which $L$ is a vector space, $\alpha$ is a linear self-map of $L$, and the skew-symmetric bilinear bracket satisfies an $\alpha$-twisted variant of the Jacobi identity, called the \emph{Hom-Jacobi identity} \eqref{eq:HomLieaxioms}.  Lie algebras are special cases of Hom-Lie algebras in which $\alpha$ is the identity map.  Some $q$-deformations of the Witt and the Virasoro algebras have the structure of a Hom-Lie algebra \cite{hls}.  Hom-Lie algebras are closely related to discrete and deformed vector fields and differential calculus \cite{hls,ls,ls2}.

Hom-Lie algebras are also useful in mathematical physics.  In \cite{yau5,yau6}, applications of Hom-Lie algebras to a generalization of the Yang-Baxter equation (YBE) \cite{baxter1,baxter2,yang} and to braid group representations \cite{artin1,artin2} are discussed.  In particular, for a vector space $M$ and a linear self-map $\alpha$ on $M$, the \emph{Hom-Yang-Baxter equation} (HYBE) for $(M,\alpha)$ is the equation
\[
(\alpha \otimes B) \circ (B \otimes \alpha) \circ (\alpha \otimes B) = (B \otimes \alpha) \circ (\alpha \otimes B) \circ (B \otimes \alpha),
\]
where $B \colon M^{\otimes 2} \to M^{\otimes 2}$ is a bilinear map that commutes with $\alpha^{\otimes 2}$.  The YBE is the special case of the HYBE with $\alpha = Id$.  It is proved in \cite[Theorem 1.1]{yau5} that every Hom-Lie algebra $(L,[-,-],\alpha)$ gives rise to a solution $B_\alpha$ of the HYBE for $(\mathbb{K} \oplus L, Id \oplus \alpha)$.  More precisely, this operator $B_\alpha$ is defined as
\[
B_\alpha((a,x) \otimes (b,y)) = (b,\alpha(y)) \otimes (a,\alpha(x)) + (1,0) \otimes (0,[x,y])
\]
for $a,b \in \mathbb{K}$ (the characteristic $0$ ground field) and $x,y \in L$. If in addition $\alpha$ is invertible, then so is $B_\alpha$.   It is also shown in \cite[Theorem 1.4]{yau5} that every solution of the HYBE for $(M,\alpha)$ gives rise to operators $B_i$ $(1 \leq i \leq n-1)$ on $M^{\otimes n}$ that satisfy the braid relations.  In particular, if $\alpha$ is invertible in the Hom-Lie algebra $(L,[-,-],\alpha)$, then these operators $B_i$ on $(\mathbb{K} \oplus L)^{\otimes n}$ satisfy the braid relations and are invertible.  So we obtain a representation of the braid group on $n$ strands on the linear automorphism group of $(\mathbb{K} \oplus L)^{\otimes n}$.  It is, therefore, useful to have concrete examples of Hom-Lie algebras.

Meanwhile in \cite{ms}, Makhlouf and Silvestrov introduced the notion of a \emph{Hom-associative algebra} $(A, \mu, \alpha)$, in which $\alpha$ is a linear self-map of the vector space $A$ and the bilinear operation $\mu$ satisfies an $\alpha$-twisted version of associativity \eqref{eq:HAaxiom}.  Associative algebras are special cases of Hom-associative algebras in which $\alpha$ is the identity map.  A Hom-associative algebra $A$ gives rise to a Hom-Lie algebra $HLie(A)$ via the commutator bracket \cite{ms}.  In this sense, Hom-associative algebras play the role of associative algebras in the Hom-Lie setting.  Moreover, any Hom-Lie algebra $L$ has a corresponding enveloping Hom-associative algebra $\UHLie(L)$ in such a way that $HLie$ and $\UHLie$ are adjoint functors \cite{yau}.  In fact, a unital version of $\UHLie(L)$ has the structure of a Hom-bialgebra \cite[Theorem 3.12]{yau2}, generalizing the bialgebra structure on the usual enveloping algebra of a Lie algebra.  Besides \cite{ms,yau,yau2,yau5,yau6}, Hom-algebras have been further studied in \cite{am,atm,cg,fg,fg2,gohr,larsson,ms2,ms3,ms4,yau4,yau7,yau8,yau9}.

There are two main purposes of this paper:
\begin{enumerate}
\item
We show how certain Hom-algebras arise naturally from classical algebras.  In particular, we show how arbitrary associative and Lie algebras deform into Hom-associative and Hom-Lie algebras, respectively, via any algebra endomorphisms.  This construction actually applies more generally to $G$-Hom-associative algebras (Theorem ~\ref{thm:algdef}), which are introduced in \cite{ms}.  This gives a systematic method for constructing many different types of Hom-algebras, including Hom-Lie algebras.
\item
We lay the foundation of a homology theory for Hom-Lie algebras.  In particular, we construct a Chevalley-Eilenberg type homology theory for Hom-Lie algebras with non-trivial coefficients.  When applied to a Lie algebra $L$, our homology of $L$ coincides with the usual Chevalley-Eilenberg homology of $L$ \cite{ce}.  The corresponding cohomology theory for Hom-Lie algebras was studied in \cite{ms3}.
\end{enumerate}

\subsection{Organization}

The rest of this paper is organized as follows.

In the next section, basic definitions about $G$-Hom-associative algebras are recalled.  It is then shown that $G$-associative algebras deform into $G$-Hom-associative algebras via an algebra endomorphism (Theorem \ref{thm:algdef}).  The desired deformations of associative and Lie algebras into their Hom counterparts are special cases of this result (Corollary ~\ref{cor:G}).  Examples of such Hom-associative and Hom-Lie deformations are then given (Examples \ref{ex:poly} - \ref{ex:Witt}).  Note that, since the appearance of an earlier version of this paper \cite{yau3}, Theorem \ref{thm:algdef} has been applied and generalized in \cite[Theorem 2.7]{am}, \cite[Theorems 1.7 and 2.6]{atm}, \cite[Section 2]{fg2}, \cite[Proposition 1]{gohr}, \cite[Theorem 3.15 and Proposition 3.30]{ms4}, \cite[Example 3.7 and Proposition 4.2]{yau2}, and \cite{yau4}--\cite{yau9}.

In Section \ref{sec:HomLie}, the homology of a Hom-Lie algebra with non-trivial coefficients is constructed (section \ref{subsec:CE}).  An interpretation of the $0$th homology module is given \eqref{eq:H0Lie}.


\section{$G$-Hom-associative algebras as deformations of $G$-associative algebras}
\label{sec:prelim}

The purposes of this section are to recall some basic definitions about $G$-Hom-associative algebras and to show that $G$-associative algebras deform into $G$-Hom-associative algebras via algebra endomorphisms (Theorem \ref{thm:algdef} and Examples \ref{ex:poly} - \ref{ex:Witt}).

\subsection{Conventions}

Throughout the rest of this paper, we work over a fixed field $\mathbb{K}$ of characteristic $0$.  Tensor products, $\Hom$, modules, and chain complexes are all meant over $\mathbb{K}$, unless otherwise specified.

\subsection{Hom-modules}
\label{subsec:Hmodules}

A \emph{Hom-module} is a pair $(M,\alpha_M)$ consisting of (i) a vector space $M$ and
(ii) a linear self-map $\alpha_M \colon M \to M$.  A \emph{morphism} $f \colon (M,\alpha_M) \to (N,\alpha_N)$ of Hom-modules is a linear map $f \colon M \to N$ such that $f \circ \alpha_M = \alpha_N \circ f$.

\subsection{$G$-Hom-associative algebras}
\label{subsec:HA}

Let $G$ be a subgroup of $\Sigma_3$, the symmetric group on three letters.  A \textbf{$G$-Hom-associative algebra} \cite{ms} is a triple $(A,\mu,\alpha)$ in which $A$ is a vector space, $\mu \colon A^{\otimes 2} \to A$ is a bilinear map, and $\alpha \colon A \to A$ is a linear map, satisfying the following \emph{$G$-Hom-associativity} axiom:
\begin{equation}
\label{eq:GHom}
\sum_{\sigma \in G} (-1)^{\varepsilon(\sigma)} \left\{(x_{\sigma(1)}x_{\sigma(2)})\alpha(x_{\sigma(3)}) - \alpha(x_{\sigma(1)})(x_{\sigma(2)}x_{\sigma(3)})\right\} = 0
\end{equation}
for $x_i \in A$, where $\varepsilon(\sigma)$ is the signature of $\sigma$.  A \textbf{$G$-associative algebra} is a not-necessarily associative algebra $(A,\mu)$, satisfying \eqref{eq:GHom} with $\alpha = Id$.  Here and in what follows, we use the abbreviation $xy$ for $\mu_A(x,y)$.  Note that in \cite{ms}, $\alpha$ was not required to be multiplicative.  In some statements in this article the multiplicativity of $\alpha$ is essential.  In all such cases,
this will be explicitly indicated.

A \emph{morphism} $f \colon (A,\mu_A,\alpha_A) \to (B,\mu_B,\alpha_B)$ of $G$-Hom-associative algebras is a morphism $f \colon (A,\alpha_A) \to (B,\alpha_B)$ of Hom-modules such that $f \circ \mu_A = \mu_B \circ f^{\otimes 2}$.

Special cases of $G$-Hom-associative algebras include the following:
\begin{enumerate}
\item
A \textbf{Hom-associative algebra} is a $G$-Hom-associative algebra in which $G$ is the trivial subgroup $\{e\}$.  The $G$-Hom-associativity axiom \eqref{eq:GHom} now takes the form
   \begin{equation}
   \label{eq:HAaxiom}
   (xy)\alpha(z) = \alpha(x)(yz),
   \end{equation}
which we call \emph{Hom-associativity}.
\item
A \textbf{Hom-Lie algebra} is a $G$-Hom-associative algebra $(A,\mu,\alpha)$ in which $\mu = [-,-]$ is skew-symmetric and $G$ is the three-element subgroup $A_3$ of $\Sigma_3$.  The $A_3$-Hom-associativity axiom \eqref{eq:GHom} is equivalent to
   \begin{equation}
   \label{eq:HomLieaxioms}
   \lbrack \alpha(x), \lbrack y, z \rbrack\rbrack + \lbrack \alpha(z), \lbrack x, y \rbrack \rbrack + \lbrack \alpha(y), \lbrack z, x \rbrack\rbrack = 0,
   \end{equation}
called the \emph{Hom-Jacobi identity}.
\item
A \textbf{Hom-left-symmetric algebra} is a $G$-Hom-associative algebra in which $G = \{e,(1~2)\}$.  The $\{e,(1~2)\}$-Hom-associativity axiom \eqref{eq:GHom} is equivalent to
\begin{equation}
\label{eq:hls}
(xy)\alpha(z) - \alpha(x)(yz) = (yx)\alpha(z) - \alpha(y)(xz).
\end{equation}
Left-symmetric algebras (also called left pre-Lie algebras and Vinberg algebras) are exactly the $\{e,(1~2)\}$-associative algebras.  In other words, left-symmetric algebras are the algebras that satisfy \eqref{eq:hls} with $\alpha = Id$.
\item
A \textbf{Hom-Lie-admissible algebras} is a $\Sigma_3$-Hom-associativity algebra.  Every $G$-Hom-associative algebra is also a Hom-Lie-admissible algebra.  Moreover, if $(A,\mu,\alpha)$ is a Hom-Lie-admissible algebra, then $(A,[-,-],\alpha)$ is a Hom-Lie algebra \cite[Section 2]{ms}, where $[-,-]$ is the commutator bracket defined by $\mu$.  A Lie-admissible algebra is exactly a $\Sigma_3$-associative algebra, i.e., a Hom-Lie-admissible algebra in which $\alpha = Id$.  Equivalently, a Lie-admissible algebra is an algebra whose commutator bracket satisfies the Jacobi identity.
\end{enumerate}

The following result says that $G$-associative algebras deform into $G$-Hom-associative algebras along any algebra endomorphism.

\begin{thm}
\label{thm:algdef}
Let $(A,\mu)$ be a $G$-associative algebra and $\alpha \colon A \to A$ be a linear map such that $\alpha \circ \mu = \mu \circ \alpha^{\otimes 2}$.  Then $(A,\mu_\alpha = \alpha \circ \mu,\alpha)$ is a $G$-Hom-associative algebra.  Moreover, $\alpha$ is multiplicative with respect to $\mu_\alpha$, i.e., $\alpha \circ \mu_\alpha = \mu_\alpha \circ \alpha^{\otimes 2}$.

Suppose that $(B, \mu')$ is another $G$-associative algebra and that $\alpha' \colon B \to B$ is a linear map such that $\alpha' \circ \mu' = \mu' \circ \alpha'^{\otimes 2}$.  If $f \colon A \to B$ is an algebra morphism (i.e., $f \circ \mu = \mu' \circ f^{\otimes 2}$) that satisfies $f \circ \alpha = \alpha' \circ f$, then $f \colon (A, \mu_\alpha, \alpha) \to (B, \mu'_{\alpha'} = \alpha' \circ \mu', \alpha')$ is a morphism of $G$-Hom-associative algebras.
\end{thm}

We will use the following observations in the proof of Theorem ~\ref{thm:algdef}.

\begin{lemma}
\label{lem:mu}
Let $A = (A,\mu)$ be a not-necessarily associative algebra and $\alpha \colon A \to A$ be an algebra morphism.  Then the multiplication $\mu_\alpha = \alpha \circ \mu$ satisfies
\begin{equation}
\label{eq:mualpha}
\mu_\alpha(\mu_\alpha(x,y),\alpha(z)) = \alpha^2((xy)z) \quad\text{and}\quad
\mu_\alpha(\alpha(x),\mu_\alpha(y,z)) = \alpha^2(x(yz))
\end{equation}
for $x,y,z \in A$, where $\alpha^2 = \alpha \circ \alpha$.  Moreover, $\alpha$ is multiplicative with respect to $\mu_\alpha$, i.e., $\alpha \circ \mu_\alpha = \mu_\alpha \circ \alpha^{\otimes 2}$.
\end{lemma}

\begin{proof}
Using the hypothesis that $\alpha$ is an algebra morphism, we have
\[
\mu_\alpha(\mu_\alpha(x,y),\alpha(z)) = \alpha(\alpha(xy)\alpha(z)) = \alpha^2((xy)z),
\]
proving the first assertion in \eqref{eq:mualpha}.  The other assertion in \eqref{eq:mualpha} is proved similarly.  For the last assertion, observe that both $\alpha \circ \mu_\alpha$ and $\mu_\alpha \circ \alpha^{\otimes 2}$ are equal to $\alpha \circ \mu \circ \alpha^{\otimes 2}$.
\end{proof}

\begin{proof}[Proof of Theorem ~\ref{thm:algdef}]
By Lemma ~\ref{lem:mu}, $\alpha$ is multiplicative with respect to $\mu_\alpha$.  Next we check \eqref{eq:GHom} with the multiplication $\mu_\alpha = \alpha \circ \mu$.  We compute as follows:
\[
\begin{split}
\sum_{\sigma \in G} & (-1)^{\varepsilon(\sigma)} \left\{\mu_\alpha(\mu_\alpha(x_{\sigma(1)},x_{\sigma(2)}),\alpha(x_{\sigma(3)})) - \mu_\alpha(\alpha(x_{\sigma(1)}),\mu_\alpha(x_{\sigma(2)},x_{\sigma(3)}))\right\}\\
& = \alpha^2\left\{\sum_{\sigma \in G} (-1)^{\varepsilon(\sigma)} \left\{(x_{\sigma(1)}x_{\sigma(2)})x_{\sigma(3)} - x_{\sigma(1)}(x_{\sigma(2)}x_{\sigma(3)})\right\}\right\} = 0.
\end{split}
\]
The first equality follows from Lemma ~\ref{lem:mu} and the linearity of $\alpha$.  The second equality follows from the hypothesis that $(A,\mu)$ is $G$-associative.

Finally, $f$ is a morphism of $G$-Hom-associative algebras because $f \circ \alpha = \alpha' \circ f$ by hypothesis and
\[
f \circ \mu_\alpha = f \circ \alpha \circ \mu = \alpha' \circ f \circ \mu = \mu'_{\alpha'} \circ f^{\otimes 2}
\]
by the assumption that $f$ is an algebra morphism.
\end{proof}

If we take $G$ to be the subgroups $\{e\}$, $A_3$, $\{e,(1~2)\}$, and $\Sigma_3$, respectively, in Theorem ~\ref{thm:algdef}, we obtain the following result.

\begin{corollary}
\label{cor:G}
Let $A = (A,\mu)$ be a not-necessarily associative algebra and $\alpha \colon A \to A$ be an algebra morphism.  Write $A_\alpha$ for the triple $(A,\mu_\alpha = \alpha \circ \mu,\alpha)$.
\begin{enumerate}
\item
If $A$ is an associative algebra, then $A_\alpha$ is a Hom-associative algebra.
\item
If $A$ is a Lie algebra, then $A_\alpha$ is a Hom-Lie algebra.
\item
If $A$ is a left-symmetric algebra, then $A_\alpha$ is a Hom-left-symmetric algebra.
\item
If $A$ is a Lie-admissible algebra, then $A_\alpha$ is a Hom-Lie-admissible algebra.
\end{enumerate}
\end{corollary}

In view of Theorem \ref{thm:algdef}, we think of the $G$-Hom-associative algebra $A_\alpha = (A, \mu_\alpha, \alpha)$ as a deformation of the $G$-associative algebra $A$ that reduces to $A$ when $\alpha = \Id_A$.  In the rest of this section, we give several examples of this kind of Hom-associative and Hom-Lie deformations.

\begin{example}[\textbf{Polynomial Hom-associative algebras}]
\label{ex:poly}
Consider the polynomial algebra $A = \mathbb{K} \lbrack x_1, \ldots , x_n \rbrack$ in $n$ variables.  Then an algebra endomorphism $\alpha$ of $A$ is uniquely determined by the $n$ polynomials $\alpha(x_i) = \sum \lambda_{i; r_1, \ldots , r_n} x_1^{r_1} \cdots x_n^{r_n}$ for $1 \leq i \leq n$.  Define $\mu_\alpha$ by
   \[
   \mu_\alpha(f,g) = f(\alpha(x_1), \ldots , \alpha(x_n)) g(\alpha(x_1), \ldots , \alpha(x_n))
   \]
for $f$ and $g$ in $A$.   By Corollary ~\ref{cor:G}, $A_\alpha = (A, \mu_\alpha, \alpha)$ is a Hom-associative algebra that reduces to the original polynomial algebra $A$ when $\alpha(x_i) = x_i$
for $1 \leq i \leq n$, i.e., $\alpha = \Id$.  We think of the collection $\left \lbrace A_\alpha \colon \alpha \text{ an algebra endomorphism of } A \right\rbrace$ as a family of deformations of the polynomial algebra $A$ into Hom-associative algebras.  A generalization of this example is considered in \cite[Example 3.32]{ms4}.\qed
\end{example}

\begin{example}[\textbf{Group Hom-associative algebras}]
Let $A = \mathbb{K}\lbrack G \rbrack$ be the group-algebra over a group $G$.  If $\alpha \colon G \to G$ is a group morphism, then it can be extended to an algebra endomorphism of $A$ by setting $\alpha \left(\sum_{g \in G} a_g g\right) = \sum_{g \in G} a_g \alpha(g)$.  By Corollary ~\ref{cor:G}, $A_\alpha = (A, \mu_\alpha, \alpha)$ is a Hom-associative algebra in which
\[
\mu_\alpha\left(\sum a_g g, \sum b_g g\right) = \sum c_g \alpha(g),
\]
where $(\sum a_g g)(\sum b_g g) = \sum c_g g$.  We think of the collection $\left\lbrace A_\alpha \colon \alpha \colon G \to G \text{ a group morphism}\right\rbrace$ as a $1$-parameter family of deformations of the group-algebra $A$ into Hom-associative algebras.  A generalization of this example is considered in \cite[Example 3.31]{ms4}.\qed
\end{example}

\begin{example}[\textbf{Hom-associative deformations by inner automorphisms}]
\label{ex:inneraut}
Let $A$ be a unital associative algebra.  Suppose that $u \in A$ is an invertible element.  Then the map $\alpha(u) \colon A \to A$ defined by $\alpha(u)(x) = uxu^{-1}$ for $x \in A$ is an algebra automorphism.  In this case, we have
   \[
   \mu_{\alpha(u)}(x,y) = uxyu^{-1}
   \]
for $x, y \in A$.  By Corollary ~\ref{cor:G}, the triple $A_u = (A, \mu_{\alpha(u)}, \alpha(u))$ is a Hom-associative algebra.  We think of the collection $\left \lbrace A_u \colon u \in A \text{ invertible} \right\rbrace$ as a $1$-parameter family of deformations of $A$ into Hom-associative algebras.\qed
\end{example}

\begin{example}[\textbf{Hom-associative deformations by nilpotent derivations}]
\label{ex:assexp}
Let $A$ be an associative algebra.  Recall that a \emph{derivation on $A$} is a linear self-map $D$ on $A$ that satisfies the Leibniz identity, $D(xy) = D(x)y + xD(y)$, for $x, y \in A$.  Such a derivation is said to be \emph{nilpotent} if $D^n = 0$ for some $n \geq 1$.  For example, if $x \in A$ is a nilpotent element, say, $x^n = 0$, then the linear self-map $\ad(x)$ on $A$ defined by $\ad(x)(y) = xy - yx$ is a nilpotent derivation on $A$.  Given a nilpotent derivation $D$ on $A$ (with, say, $D^n = 0$), the linear self-map
\[
\exp D = \Id_A + D + \frac{1}{2}D^2 + \cdots + \frac{1}{(n-1)!}D^{n-1}
\]
is actually an algebra automorphism of $A$ (see, e.g., \cite[p.26]{abe}).  With $\mu_{\exp D}$ defined as $\mu_{\exp D}(x,y) = (\exp D)(xy)$, Corollary ~\ref{cor:G} shows that we have a Hom-associative algebra $A_D = (A, \mu_{\exp D}, \exp D)$.  We think of the collection $\left\lbrace A_D \colon D \text{ is a nilpotent derivation on $A$}\right\rbrace$ as a $1$-parameter family of deformations of $A$ into Hom-associative algebras.\qed
\end{example}

\begin{example}[\textbf{Hom-Lie $\mathsf{sl(2,\bbC)}$}]
\label{ex:sl2}
Consider the complex Lie algebra $\mathsf{sl(2,\bbC)}$ of $2 \times 2$ matrices with trace $0$.  A standard linear basis of $\mathsf{sl(2,\bbC)}$ consists of the elements
   \[
   h = \begin{pmatrix} 1 & 0 \\ 0 & -1 \end{pmatrix}, \quad
   e = \begin{pmatrix} 0 & 1 \\ 0 & 0 \end{pmatrix}, \quad \text{and} \quad
   f = \begin{pmatrix} 0 & 0 \\ 1 & 0 \end{pmatrix},
   \]
which satisfy the relations $\lbrack h,e \rbrack = 2e$, $\lbrack h,f \rbrack = -2f$, and $\lbrack e,f \rbrack = h$.  Let $\lambda \not= 0$ be a scalar in $\bbC$.  Consider the linear map $\alpha_\lambda \colon \mathsf{sl(2,\bbC)} \to \mathsf{sl(2,\bbC)}$ defined by
   \[
   \alpha_\lambda(h) = h, \quad
   \alpha_\lambda(e) = \lambda e, \quad \text{and} \quad
   \alpha_\lambda(f) = \lambda^{-1}f
   \]
on the basis elements.  The map $\alpha_\lambda$ is actually a Lie algebra morphism.  In fact, it suffices to check this on the basis elements, which is immediate from the definition of $\alpha_\lambda$.  By Corollary ~\ref{cor:G}, we have a Hom-Lie algebra $\mathsf{sl(2,\bbC)}_\lambda = (\mathsf{sl(2,\bbC)}, \lbrack -,- \rbrack_{\alpha_\lambda}, \alpha_\lambda)$.  The Hom-Lie algebra bracket $\lbrack -,- \rbrack_{\alpha_\lambda}$ on the basis elements is given by
   \[
   \lbrack h,e \rbrack_{\alpha_\lambda} = 2\lambda e, \quad
   \lbrack h,f \rbrack_{\alpha_\lambda} = -2\lambda^{-1}f, \quad \text{and} \quad
   \lbrack e,f \rbrack_{\alpha_\lambda} = h.
   \]
We think of the collection $\left\lbrace \mathsf{sl(2,\bbC)}_\lambda \colon \lambda \not=0 \text{ in $\bbC$} \right\rbrace$ as a one-parameter family of deformations of $\mathsf{sl(2,\bbC)}$ into Hom-Lie algebras.\qed
\end{example}

\begin{example}[\textbf{Hom-Lie $\mathsf{sl(n,\bbC)}$}]
\label{ex:sln}
This is a generalization of the previous example to $n > 2$.  Let $\mathsf{sl(n,\bbC)}$ be the complex Lie algebra of $n \times n$ matrices with trace $0$.  It is generated as a Lie algebra by the elements
   \[
   e_i = E_{i,i+1}, \quad
   f_i = E_{i+1,i}, \quad \text{and} \quad
   h_i = E_{ii} - E_{i+1,i+1}
   \]
for $1 \leq i \leq n-1$, where $E_{ij}$ denotes the matrix with $1$ in the $(i,j)$-entry and $0$ everywhere else.  These elements satisfy some relations similar to those of $\mathsf{sl(2,\bbC)}$ (see, e.g., \cite[p.9]{kang}).

Let $\lambda_1, \ldots, \lambda_{n-1}$ be non-zero scalars in $\bbC$.  Consider the map $\alpha_{\lambda_1, \ldots , \lambda_{n-1}} \colon \mathsf{sl(n,\bbC)} \to \mathsf{sl(n,\bbC)}$ defined on the generators by
   \[
   \alpha_{\lambda_1, \ldots , \lambda_{n-1}}(e_i) = \lambda_i e_i, \quad
   \alpha_{\lambda_1, \ldots , \lambda_{n-1}}(f_i) = \lambda_i^{-1} f_i, \quad \text{and} \quad
   \alpha_{\lambda_1, \ldots , \lambda_{n-1}}(h_i) = h_i
   \]
for $1 \leq i \leq n-1$.  It is easy to check that $\alpha_{\lambda_1, \ldots , \lambda_{n-1}}$ actually defines a Lie algebra morphism.  By Corollary ~\ref{cor:G}, we have a Hom-Lie algebra
   \[
   \mathsf{sl(n,\bbC)}_{\lambda_1, \ldots , \lambda_{n-1}}
   = (\mathsf{sl(n,\bbC)}, \lbrack -,- \rbrack_{\alpha_{\lambda_1, \ldots , \lambda_{n-1}}}, \alpha_{\lambda_1, \ldots , \lambda_{n-1}}).
   \]
We think of the collection $\left \lbrace \mathsf{sl(n,\bbC)}_{\lambda_1, \ldots , \lambda_{n-1}} \colon \lambda_1, \ldots , \lambda_{n-1} \not= 0 \text{ in $\bbC$} \right\rbrace$ as an $(n-1)$-parameter family of deformations of $\mathsf{sl(n,\bbC)}$ into Hom-Lie algebras.\qed
\end{example}

\begin{example}[\textbf{Hom-Lie Heisenberg algebra}]
\label{ex:Heisenberg}
Let $\mathsf{H}$ be the $3$-dimensional Heisenberg Lie algebra, which consists of the strictly upper-triangular complex $3 \times 3$ matrices.  It has a standard linear basis consisting of the elements
   \[
   e = \begin{pmatrix} 0 & 1 & 0 \\ 0 & 0 & 0 \\ 0 & 0 & 0 \end{pmatrix}, \quad
   f = \begin{pmatrix} 0 & 0 & 0 \\ 0 & 0 & 1 \\ 0 & 0 & 0 \end{pmatrix}, \quad \text{and} \quad
   h = \begin{pmatrix} 0 & 0 & 1 \\ 0 & 0 & 0 \\ 0 & 0 & 0 \end{pmatrix}.
   \]
The Heisenberg relation $\lbrack e,f \rbrack = h$ is satisfied, and $\lbrack e,h \rbrack = 0 = \lbrack f,h \rbrack$ are the other two relations for the basis elements.

Let $\lambda_1$ and $\lambda_2$ be non-zero scalars in $\bbC$.  Consider the map $\alpha_{\lambda_1,\lambda_2} \colon \mathsf{H} \to \mathsf{H}$
defined on the basis elements by
   \[
   \alpha_{\lambda_1,\lambda_2}(e) = \lambda_1 e, \quad
   \alpha_{\lambda_1,\lambda_2}(f) = \lambda_2 f, \quad \text{and} \quad
   \alpha_{\lambda_1,\lambda_2}(h) = \lambda_1\lambda_2 h.
   \]
It is straightforward to check that $\alpha_{\lambda_1,\lambda_2}$ defines a Lie algebra morphism.  By Corollary ~\ref{cor:G}, we have a Hom-Lie algebra $\mathsf{H}_{\lambda_1,\lambda_2} = (\mathsf{H}, \lbrack -,- \rbrack_{\alpha_{\lambda_1,\lambda_2}}, \alpha_{\lambda_1,\lambda_2})$,
whose bracket satisfies the \emph{twisted Heisenberg relation}
   \[
   \lbrack e,f \rbrack_{\alpha_{\lambda_1,\lambda_2}} = \lambda_1\lambda_2 h.
   \]
We think of the collection $\left\lbrace \mathsf{H}_{\lambda_1,\lambda_2} \colon \lambda_1, \lambda_2 \in \bbC \right\rbrace$ as a $2$-parameter family of deformations of $\mathsf{H}$ into Hom-Lie algebras.\qed
\end{example}

\begin{example}[\textbf{Matrix Hom-Lie algebras}]
\label{ex:matrix}
Let $G$ be a matrix Lie group (e.g., $\mathsf{GL(n,\bbC)}$, $\mathsf{SL(n,\bbC)}$, $\mathsf{U(n)}$, $\mathsf{O(n)}$, and $\mathsf{Sp(n)}$), and let $\mathfrak{g}$ be the Lie algebra of $G$.  Given any element $x \in G$, it is well-known that the map $Ad_x \colon \mathfrak{g} \to \mathfrak{g}$ defined by $Ad_x(g) = xgx^{-1}$ is a Lie algebra morphism (see, e.g., \cite[Proposition 2.23]{hall}).  By Corollary ~\ref{cor:G}, we have a Hom-Lie algebra
   \[
   \mathfrak{g}_x = (\mathfrak{g}, \lbrack -,- \rbrack_{Ad_x}, Ad_x)
   \]
in which
   \[
   \lbrack g_1, g_2 \rbrack_{Ad_x} = x(g_1g_2 - g_2g_1)x^{-1}
   \]
for $g_1, g_2 \in \mathfrak{g}$.  We think of the collection $\left\lbrace \mathfrak{g}_x \colon x \in G \right\rbrace$ as a $1$-parameter family of deformations of $\mathfrak{g}$ into Hom-Lie algebras.\qed
\end{example}

\begin{example}[\textbf{Hom-Lie Witt algebra}]
\label{ex:Witt}
The Witt algebra $W$ is the complex Lie algebra of derivations on the Laurent polynomial algebra $\bbC \lbrack t^{\pm 1} \rbrack$.  It can be regarded as the one-dimensional $\bbC \lbrack t^{\pm 1} \rbrack$-module
   \[
   W = \bbC \lbrack t^{\pm 1} \rbrack \cdot \frac{d}{dt},
   \]
whose Lie bracket is given by
   \[
   \left\lbrack f \cdot \frac{d}{dt}, g \cdot \frac{d}{dt} \right\rbrack
   = \left(f \frac{dg}{dt} - g\frac{df}{dt}\right)\cdot \frac{d}{dt}
   \]
for $f, g \in \bbC \lbrack t^{\pm 1} \rbrack$ (see, e.g., \cite[Example 11]{hls}).  Given any scalar $\lambda \in \bbC$, the map $\alpha_\lambda \colon W \to W$ defined by
   \[
   \alpha_\lambda\left(f \cdot \frac{d}{dt}\right) = f(\lambda + t) \cdot \frac{d}{dt}
   \]
is easily seen to be a Lie algebra morphism.  By Corollary ~\ref{cor:G}, we have a Hom-Lie algebra $W_\lambda = \left(W, \lbrack -,- \rbrack_{\alpha_\lambda}, \alpha_\lambda\right)$,
in which the bracket is given by
   \[
   \left\lbrack f \cdot \frac{d}{dt}, g \cdot \frac{d}{dt} \right\rbrack_{\alpha_\lambda}
   = \left(f(\lambda + t)\frac{dg}{dt}(\lambda + t) - g(\lambda + t)\frac{df}{dt}(\lambda + t)\right) \cdot \frac{d}{dt}.
   \]
We think of the collection $\left\lbrace W_\lambda \colon \lambda \in \bbC \right\rbrace$ as a $1$-parameter family of deformations of the Witt algebra $W$ into Hom-Lie algebras.\qed
\end{example}


\section{Homology for Hom-Lie algebras}
\label{sec:HomLie}

The purpose of this section is to construct the homology for a Hom-Lie algebra.  We begin by defining the coefficients.

\subsection{Hom-$L$-module}
\label{subsec:HomLmodule}

From now on, $(L, \lbrack -, - \rbrack, \alpha_L)$ will denote a Hom-Lie algebra \eqref{eq:HomLieaxioms} in which $\alpha_L$ is multiplicative with respect to $[-,-]$, unless otherwise specified.

By a \emph{(right) Hom-$L$-module}, we mean a Hom-module $(M,\alpha_M)$ that comes equipped with a right $L$-action, $\rho \colon M \otimes L \to M$ $(m \otimes x \mapsto mx)$, such that the following two conditions are satisfied for $m \in M$ and $x,y \in L$:
   \begin{equation}
   \label{eq:HomLmodule}
   \begin{split}
   \alpha_M(m)\lbrack x, y \rbrack &= (mx)\alpha_L(y) - (my)\alpha_L(x), \\
   \alpha_M(mx) &= \alpha_M(m)\alpha_L(x)
   \end{split}
   \end{equation}

\begin{example}
\label{ex:HomLmodule}
Here are some examples of Hom-$L$-modules.
\begin{enumerate}
\item
One can consider $L$ itself as a Hom-$L$-module in which the $L$-action is the bracket $\lbrack -,- \rbrack$.

\item
If $\mathfrak{g}$ is a Lie algebra and $M$ is a right $\mathfrak{g}$-module in the usual sense, then $(M,\Id_M)$ is a Hom-$\mathfrak{g}$-module.\qed
\end{enumerate}
\end{example}

\subsection{The chain complex $CE^\alpha_*(L,M)$}
\label{subsec:CE}

For the rest of this section, $(M,\alpha_M)$ will denote a fixed Hom-$L$-module, where $L$ is a Hom-Lie algebra.  For $n \geq 0$, let $\Lambda^nL$ denote the $n$th exterior power of $L$, with $\Lambda^0L = \mathbb{K}$.  A typical generator in $\Lambda^nL$ is denoted by $x_1 \wedge \cdots \wedge x_n$ with each $x_i \in L$.  We will use the following abbreviations:
\[
   \begin{split}
   x_1 \cdots \widehat{x_i} \cdots x_n
   &= x_1 \wedge \cdots \wedge x_{i-1} \wedge x_{i+1} \wedge \cdots \wedge x_n, \\
   \alpha_L(x_1 \cdots \widehat{x_i} \cdots x_n)
   &= \alpha_L(x_1) \wedge \cdots \wedge \alpha_L(x_{i-1}) \wedge \alpha_L(x_{i+1}) \wedge \cdots \wedge \alpha_L(x_n).
   \end{split}
\]
Likewise, the symbols $x_1 \cdots \widehat{x_i} \cdots \widehat{x_j} \cdots x_n$, $\alpha_L(x_1 \cdots \widehat{x_i} \cdots \widehat{x_j} \cdots x_n)$, and so forth mean that the terms $\widehat{x_i}$, $\widehat{x_j}$, etc., are omitted.

Define the \emph{module of $n$-chains of $L$ with coefficients in $M$} as
   \[
   CE^\alpha_n(L,M) = M \otimes \Lambda^nL.
   \]
For $p \geq 1$, define a linear map $d_p \colon CE^\alpha_p(L,M) \to CE^{\alpha}_{p-1}(L,M)$ by setting (for $m \in M, \, x_i \in L$)
   \begin{equation}
   \label{eq:theta}
   d_p(m \otimes x_1 \wedge \cdots \wedge x_p) = \eta_1 + \eta_2,
   \end{equation}
where
   \[
   \eta_1
   = \sum_{i=1}^p (-1)^{i+1} mx_i \otimes \alpha_L(x_1 \cdots \widehat{x_i} \cdots x_p)
   \]
and
   \[
   \eta_2
   = \sum_{i<j} (-1)^{i+j} \alpha_M(m) \otimes \lbrack x_i, x_j \rbrack \wedge \alpha_L(x_1 \cdots \widehat{x_i} \cdots \widehat{x_j} \cdots x_p).
   \]

\begin{thm}
\label{thm:dCE}
The data $(CE^\alpha_*(L,M),d)$ forms a chain complex.
\end{thm}

\begin{proof}
Using the notations in \eqref{eq:theta}, we have
   \[
   \begin{split}
   d^2(m \otimes x_1 \wedge \cdots \wedge x_p)
   & = d(\eta_1) + d(\eta_2) \\
   & = (\eta_{11} + \eta_{12}) + (\eta_{21} + \eta_{22}).
   \end{split}
   \]
Therefore, to prove the Theorem, it suffices to show that
   \begin{equation}
   \label{eq:eta21}
   \eta_{11} + \eta_{12} + \eta_{21} = 0
   \end{equation}
and
   \begin{equation}
   \label{eq:eta22}
   \eta_{22} = 0.
   \end{equation}
To prove \eqref{eq:eta21}, first note that $\eta_{21}$ is a sum of $p-1$ terms, the first of which is
   \begin{equation}
   \label{eq:eta21i=1}
   \sum_{i<j} (-1)^{i+j}\left(\alpha_M(m)\lbrack x_i,x_j \rbrack\right) \otimes \alpha_L^2(x_1 \cdots \widehat{x_i} \cdots \widehat{x_j} \cdots x_p).
   \end{equation}
On the other hand, we have that
   \[
   \begin{split}
   \eta_{11}
   & = \sum_{i=1}^p (-1)^{i+1} \sum_{j < i} (-1)^{j+1} (mx_i)\alpha_L(x_j) \otimes \alpha_L^2(x_1 \cdots \widehat{x_j} \cdots \widehat{x_i} \cdots x_p) \\
   & \relphantom{} + \sum_{i=1}^p (-1)^{i+1} \sum_{j>i} (-1)^{(j-1)+1} (mx_i)\alpha_L(x_j) \otimes \alpha_L^2(x_1 \cdots \widehat{x_i} \cdots \widehat{x_j} \cdots x_p) \\
   & = - \sum_{i<j} (-1)^{i+j} \left((mx_i)\alpha_L(x_j) - (mx_j)\alpha_L(x_i)\right) \otimes \alpha_L^2(x_1 \cdots \widehat{x_i} \cdots \widehat{x_j} \cdots x_p).
   \end{split}
   \]
By the first Hom-$L$-module axiom \eqref{eq:HomLmodule}, the last line is equal to \eqref{eq:eta21i=1} with a minus sign.

The other $p-2$ terms in $\eta_{21}$ are given by the sum
   \begin{multline}
   \label{eq:eta21i>1}
   \sum_{i<j<k} (-1)^{i+j+k}
   (\alpha_M(m) \alpha_L(x_i) \otimes \alpha_L(\lbrack x_j,x_k \rbrack) \wedge z
   - \alpha_M(m)\alpha_L(x_j) \otimes \alpha_L(\lbrack x_i,x_k \rbrack) \wedge z \\
   + \alpha_M(m)\alpha_L(x_k) \otimes \alpha_L(\lbrack x_i,x_j\rbrack) \wedge z),
   \end{multline}
where
   \begin{equation}
   \label{eq:z}
   z = \alpha_L^2(x_1 \cdots \widehat{x_i} \cdots \widehat{x_j} \cdots \widehat{x_k} \cdots x_p).
   \end{equation}
Using the second Hom-$L$-module axiom \eqref{eq:HomLmodule} and the multiplicativity of $\alpha_L$, we can rewrite \eqref{eq:eta21i>1} as
   \begin{multline}
   \label{eq':eta21i>1}
   \sum_{i<j<k} (-1)^{i+j+k}
   (\alpha_M(mx_i) \otimes \lbrack \alpha_L(x_j), \alpha_L(x_k) \rbrack \wedge z
   - \alpha_M(mx_j) \otimes \lbrack \alpha_L(x_i), \alpha_L(x_k)\rbrack \wedge z \\
   + \alpha_M(mx_k) \otimes \lbrack \alpha_L(x_i), \alpha_L(x_j) \rbrack \wedge z).
   \end{multline}
It is straightforward to see that \eqref{eq':eta21i>1} is equal to $\eta_{12}$ with a minus sign.  So far we have proved \eqref{eq:eta21}.

To show \eqref{eq:eta22}, observe that
   \begin{equation}
   \label{eq':eta22}
   \eta_{22} = \sum_{i<j<k} (-1)^{i+j+k} \alpha_M^2(m) \otimes y \wedge z
   + \sum_{i<j<k<l} (-1)^{i+j+k+l} \alpha_M^2(m) \otimes u \wedge w,
   \end{equation}
where $z$ is as in \eqref{eq:z} and
   \[
   \begin{split}
   y
   & = \lbrack \lbrack x_i,x_j \rbrack, \alpha_L(x_k) \rbrack
   + \lbrack \lbrack x_j,x_k \rbrack, \alpha_L(x_i) \rbrack
   + \lbrack \lbrack x_k,x_i \rbrack, \alpha_L(x_j) \rbrack, \\
   u
   & = \lbrack \alpha_L(x_i), \alpha_L(x_j)\rbrack \wedge \alpha_L(\lbrack x_k,x_l \rbrack)
   + \lbrack \alpha_L(x_k), \alpha_L(x_l) \rbrack \wedge \alpha_L(\lbrack x_i,x_j \rbrack) \\
   & \relphantom{} - \lbrack \alpha_L(x_i), \alpha_L(x_k) \rbrack \wedge \alpha_L(\lbrack x_j,x_l\rbrack)
   - \lbrack \alpha_L(x_j), \alpha_L(x_l) \rbrack \wedge \alpha_L(\lbrack x_i,x_k\rbrack) \\
   & \relphantom{} + \lbrack \alpha_L(x_i), \alpha_L(x_l) \rbrack \wedge \alpha_L(\lbrack x_j,x_k \rbrack)
   + \lbrack \alpha_L(x_j), \alpha_L(x_k) \rbrack \wedge \alpha_L(\lbrack x_i,x_l\rbrack), \\
   w
   & = \alpha_L^2(x_1 \cdots \widehat{x_i} \cdots \widehat{x_j} \cdots \widehat{x_k} \cdots \widehat{x_l} \cdots x_p).
   \end{split}
   \]
It follows from the Hom-Jacobi identity \eqref{eq:HomLieaxioms} and the skew-symmetry of $[-,-]$ that
   \begin{equation}
   \label{eq:y=0}
   y = 0.
   \end{equation}
Likewise, using the multiplicativity of $\alpha_L$ and that $a \wedge b = - b \wedge a$ in an exterior algebra, one infers that
   \begin{equation}
   \label{eq:u=0}
   u = 0.
   \end{equation}
Combining \eqref{eq':eta22}, \eqref{eq:y=0}, and \eqref{eq:u=0}, it follows that $\eta_{22} = 0$, which proves \eqref{eq:eta22}.
\end{proof}

\subsection{Homology}
\label{subsec:homologyHLie}

In view of Theorem \ref{thm:dCE}, we define the \emph{$n$th homology of $L$ with coefficients in $M$} as
   \[
   H^\alpha_n(L,M) = H_n(CE^\alpha_*(L,M)).
   \]
Note that for a Lie algebra $\mathfrak{g}$ and a right $\mathfrak{g}$-module $M$, the chain complex $CE^\alpha_*(\mathfrak{g},M)$ is exactly the Chevalley-Eilenberg complex \cite{ce} that defines the Lie algebra homology of $\mathfrak{g}$ with coefficients in the right $\mathfrak{g}$-module $M$.  This justifies our choice of notation.

\subsection{The $0$th homology module $H^\alpha_0(L,M)$}
\label{subsec:HCE0}

Since the differential $d_1 \colon M \otimes \Lambda^1L = M \otimes L \to M$ is the right $L$-action map on $M$, it follows that
   \[
   H^\alpha_0(L,M) = \frac{M}{span_{\mathbb{K}}\lbrace mx \colon m \in M,\, x \in L \rbrace}.
   \]
In particular, when $L$ is considered as a Hom-$L$-module via its bracket, we have that
   \begin{equation}
   \label{eq:H0Lie}
   H^\alpha_0(L,L) = \frac{L}{\lbrack L, L \rbrack},
   \end{equation}
which is the abelianization of $L$ with respect to its bracket.


\end{document}